\documentclass{article}
\usepackage{amssymb}
\usepackage{amsfonts}
\usepackage{amsmath}
\usepackage{chicago}

\begin{document}

\begin{center}
{\Large Some remarks about the levels and sublevels of algebras obtained by
the Cayley-Dickson process}%
\begin{equation*}
\end{equation*}%
\textbf{Cristina FLAUT}

\textquotedblright Ovidius\textquotedblright\ University of Constanta,
Romania\\[0pt]
e-mail: cflaut@univ-ovidius.ro; cristina\_flaut@yahoo.com\\[0pt]
http://cristinaflaut.wikispaces.com/%
\begin{equation*}
\end{equation*}
\end{center}

\textbf{Abstract. }\ {\small In this paper we improve the level \ and
sublevel of algebras obtained by the Cayley-Dickson process when their level
\ and sublevel are greater than dimension of the algebras.}

\bigskip

\textbf{Keywords.} Cayley--Dickson process, Division algebra. Level and
sublevel of an algebra

\textbf{Mathematics Subject Classification (2000).} 17A35, 17A20,
17A75,17A45.%
\begin{equation*}
\end{equation*}

\textbf{0. Introduction.}%
\begin{equation*}
\end{equation*}

In this paper, we assumed that the field $K$ is commutative with $charK$ $%
\neq 2$ and quadratic forms over the field $K$ are always assumed to be
finite-dimensional and nondegenerate. For the basic terminology and results
of quadratic and symmetric bilinear \ spaces, the reader is referred to
[Sch; 85].

For a given field $K,$ its \textit{level}, denoted by $s(K)$, is considered
to be the smallest natural number $n$ such that $-1$ is a sum of $n$ squares
of $K$. If $-1$ is not a sum of squares of $K$, then the level of the field $%
K$ is infinite.

In [Pf; 65], Pfister showed that if we have a finite level for a field, then
this level is a power of $2$ and any power of $\ 2$ could be realized as the
level of a field. As a generalization of this definition appears the notions
of level and sublevel of an algebra $A.$

The \textit{level} of the algebra $A$, denoted by $s(A)$, is the least
integer $n$ such that $-1$ is a sum of $n$ squares in $A$.

The \textit{sublevel} of the algebra $A$, denoted by \underline{$s$}$(A)$,
is the least integer $n$ such that $0$ is a sum of $n+1$ nonzero squares of
elements in $A$. If these numbers do not exist, then the level and sublevel
are infinite. Obviously, \underline{$s$}$(A)$ $\leq s(A)$.

There are many papers devoted to the study of the level and sublevel of
quaternion algebras, octonion algebras, composition algebras or algebras
obtained by the Cayley-Dickson process: [Hoff; 95],[Hoff; 98], [Lew; 87],
[Lew; 89], [O' Sh; 07], [O' Sh; 10], [O' Sh; 11], [Ti, Va; 87], etc. \ In
[Lew; 87], D. W. Lewis constructed quaternion division algebras of level $%
2^{k}$ and $2^{k}+1$ for all\thinspace \thinspace $k\in \mathbb{N}-\{0\}$
and he asked if there exist other values for the level of the quaternion
division algebras. As an answer of this question, in [Hoff; 08], D. W.
Hoffman proved that there are many other values, different from $2^{k}$ or $%
2^{k}+1,$ which could be realized as a level of quaternion division
algebras. He showed that \ for each $k\in \mathbb{N},$ $k\geq 2,$ there
exist quaternion division algebras $D$ with the level $s\left( D\right) \in
\lbrack 2^{k}+2,2^{k+1}-1]\cap \mathbb{N}.$

In [Pu; 05], Susanne Pumpl\"{u}n proved the existence of octonion division
algebras of level $2^{k}$ and $2^{k}+1$ for all\thinspace \thinspace $k\in 
\mathbb{N}-\{0\}$ and in [O' Sh; 10], Theorem 3.6., O'Shea constructed
octonion division algebras of level $6$ and $7.$ These values, $6$ and $7,$
are still the only known exact values for the level of octonion division
algebras, other than $2^{k}$ or $2^{k}+1,$\thinspace \thinspace $k\in 
\mathbb{N}-\{0\}.\ \ $

In the case of \ quaternion and octonion \ division algebras it is still not
known which exact numbers can be realized as their levels and sublevels. For
the integral domains, the problem of level was solved in [Da, La, Pe; 80],
when \ Z.D. Dai, \ T. Y. Lam \ and C. K. Peng, where they proved that any
positive integer $n$ could be realized as the level of an integral domain.
As a generalization of this last result, in [Fl; 13] was proved that for any
positive integer $n,$\ there is an algebra $A$\ obtained by the
Cayley-Dickson process with the norm form anisotropic over a suitable field,
which has the level $n\in N-\{0\},$ ([Fl; 13], Theorem 2.9). This result is
a better one since was replaced division property with the anisotropy of the
norm form. It is well known that for an algebra $A,$ obtained by the
Cayley-Dickson process, division property implies that its norm form is
anisotropic, but there are algebras $A~$obtained by the Cayley-Dickson
process with the norm form \ anisotropic which are not division algebras, as
for example real sedenion algebra, (see [Fl; 13], Remark 1.3, i)).

Since in the above mentioned result, the chosen algebras has the level less
or equal with their dimension, in this paper we will try to improve the
bounds of these values when the level and sublevel are greater than the
dimension of the algebra.

\begin{equation*}
\end{equation*}

\textbf{1. Preliminaries}%
\begin{equation*}
\end{equation*}

It is well known that a regular quadratic form over the field $K$ can be
diagonalized. If we \ consider $a_{1},a_{2},...a_{n}\in K^{\ast },n\in 
\mathbb{N}-\{0\},$ we denote by $<a_{1},a_{2},...a_{n}>$ the $n$-dimensional
quadratic form $a_{1}X_{1}^{2}+a_{2}X_{2}^{2}+...a_{n}X_{n}^{2}$.

A quadratic form $q:V\rightarrow K$ with the property $\,q\left( x\right) =0 
$ implies $x=0$ is called \textit{anisotropic}, otherwise $q$ is called 
\textit{isotropic}.

We denote by $p\perp q$ the orthogonal sum of the regular quadratic forms $p$
and $q$ over the field $K$ and by $p\otimes q$ their tensor product. If $%
n\in \mathbb{N}$, we will denote the orthogonal sum of $n$ copies of $q$ by $%
n\times q$.

We consider $\varphi $ a $n-$dimensional quadratic irreducible form over $K,$
$n\in N,n>1,$ which is not isometric to the hyperbolic plane, $<1,-1>.$ We
can consider $\varphi $ as a homogeneous polynomial of degree $2$,%
\begin{equation*}
\varphi \left( X\right) =\varphi \left( X_{1},...X_{n}\right) =\sum
a_{ij}X_{i}X_{j},a_{ij}\in K^{\ast }.
\end{equation*}%
The \textit{function field of }$\varphi ,\,$denoted \ by $K(\varphi ),$ is
the quotient field of the integral domain 
\begin{equation*}
K[X_{1},...,X_{n}]\,/\,\left( \varphi \left( X_{1},...,X_{n}\right) \right) .
\end{equation*}

Let $q$ be a regular quadratic form over the field $K.$ This form can be
written under the form 
\begin{equation*}
\varphi \cong \varphi _{an}\bot i_{W}\left( \varphi \right) \times
<1,-1>\cong \varphi _{an}\bot \varphi _{h},
\end{equation*}%
where the anisotropic form $q_{an}$ and the integer $i_{W}\left( \varphi
\right) $ are uniquely determined. The number $i_{W}\left( \varphi \right) $
is called the \textit{Witt index }of the quadratic form $\varphi .$ It is
clear that the Witt index of $\varphi $ is $\frac{1}{2}\dim \varphi _{h}$
and it is the dimension of a maximal totally isotropic subform of $\varphi .$
A quadratic form $q$ is hyperbolic if $q_{an}$ is trivial. \ In this
situation, the subform $\varphi _{h}$ is hyperbolic. The \textit{first Witt
index} of \ a quadratic form $\varphi $ is the Witt index of $\varphi $ over
its function field and it is denoted by $i_{1}\left( \varphi \right) .$ The 
\textit{essential dimension} of $\varphi $ is 
\begin{equation*}
\text{\textit{dim}}_{es}\left( \varphi \right) =\text{\textit{dim}}\left(
\varphi \right) -i_{1}\left( \varphi \right) +1.
\end{equation*}%
$\ \ $

Let $n\in \mathbb{N}-\{0\}$.$\,$ An $n-$\textit{fold Pfister form} over $K$
is a quadratic form of the type $\ $%
\begin{equation*}
<1,a_{1}>\otimes ...\otimes <1,a_{n}>,a_{1},...,a_{n}\in K^{\ast },
\end{equation*}%
denoted by $\ll a_{1},a_{2},...,a_{n}\gg .$ For $n\in N,n>1,\,\,\,$a Pfister
form $\varphi $ can be written under the form 
\begin{equation*}
<1,a_{1}>\otimes \text{...}\otimes <1,a_{n}>=<1,a_{1},a_{2},\text{...}%
,a_{n},a_{1}a_{2},\text{...},a_{1}a_{2}a_{3},\text{...},a_{1}a_{2}...a_{n}>.
\end{equation*}%
If $\varphi =<1>\perp \varphi ^{\prime },$ then $\varphi ^{\prime }$ is
called \textit{the pure subform} of\thinspace $\varphi .$ A Pfister form is
hyperbolic if and only if it is isotropic. Therefore a Pfister form is
isotropic if and only if its pure subform is isotropic.( See [Sch; 85] )

A quadratic form $\varphi $ is called a \textit{Pfister neighbor} if we can
find an $n$-fold Pfister form $\Phi $ such that $dim\varphi >2^{n-1}$ and $%
\alpha \Psi \simeq \varphi \perp \phi ,$ with $\phi $ a quadratic form and $%
\alpha \in K^{\ast }$

We consider the field $L$ and we define 
\begin{equation*}
L^{\infty }=L\cup \{\infty \},
\end{equation*}%
where $x+\infty =x$, for $x\in K,~x\infty =\infty $ for $x\in K^{\ast
},\infty \infty =\infty ,\frac{1}{\infty }=0,\frac{1}{0}=\infty $.

We define an $L-$\textit{place }of the field $K$ as a map $\lambda
:K\rightarrow L^{\infty }$ with the properties:%
\begin{equation*}
\lambda \left( x+y\right) =\lambda \left( x\right) +\lambda \left( y\right)
,\lambda \left( xy\right) =\lambda \left( x\right) \lambda \left( y\right) ,
\end{equation*}%
whenever the right sides are defined.

A subset $P$ of $K$ is called an \textit{ordering }of $K$ if the following
conditions are fulfilled: 
\begin{equation*}
P+P\subset P,P\cdot P\subset P,-1\notin P,
\end{equation*}%
$\ \ $%
\begin{equation*}
\{x\in K\,\,/\,x\ \ \text{is \ a sum of squares in\ \ }K\}\subset P,P\cup
-P=K,P\cap -P=0.
\end{equation*}

A field $K$ with an ordering defined on is called an \textit{ordered} field.
For $x,y\in K,$ $K$ an ordered field, we define $x>y$ \thinspace
if\thinspace $\left( x-y\right) \in P.$

If$\ \ \varphi $ $=<a_{1},...,a_{n}>$ is a quadratic form over a formally
real field $K$ and $P$ is an ordering on $K,$ \textit{the signature} of $%
\varphi $ at $P$ is 
\begin{equation*}
sgn\left( \varphi \right) =\left\vert \{i\mid ~a_{i}>_{P}0\}\right\vert
-\left\vert \{\{i\mid ~a_{i}<_{P}0\}\}\right\vert .
\end{equation*}%
The quadratic form $q$ is called \textit{indefinite} at ordering $P$ if $%
\dim \varphi >\left\vert sgn\varphi \right\vert .$

In the following, we briefly present the \textit{Cayley-Dickson process} and
some of the properties of the algebras obtained. For other details, the
reader is referred to $\left[ \text{Sc; 66}\right] $ and [Sc; 54].

Let $A$ be a finite dimensional unitary algebra over a field $\ K$ with $\,$%
\begin{equation*}
\,\,\,\overline{\phantom{x}}:A\rightarrow A,a\rightarrow \overline{a},
\end{equation*}%
$\,\,$ a linear map satisfying the following relations:$\,\,\,\,\,$%
\begin{equation*}
\overline{ab}=\overline{b}\overline{a},\,\overline{\overline{a}}=a,
\end{equation*}%
$\,\,$and 
\begin{equation*}
a+\overline{a},a\overline{a}\in K\cdot 1\ \text{for all }a,b\in A
\end{equation*}%
and called a \textit{scalar} \textit{involution. }The element $\,\overline{a}
$ is called the \textit{conjugate} of the element $a,$ the linear form$\,\,$%
\begin{equation*}
\,\,t:A\rightarrow K\,,\,\,t\left( a\right) =a+\overline{a}
\end{equation*}%
and the quadratic form 
\begin{equation*}
n:A\rightarrow K,\,\,n\left( a\right) =a\overline{a}\ 
\end{equation*}%
are called the \textit{trace} and the \textit{norm \ }of \ the element $a.$ $%
\,$

Let$\,\,\,\gamma \in K$ \thinspace be a fixed non-zero element. On the
vector space $A\oplus A,$ we define the following algebra multiplication 
\begin{equation*}
(a_{1},a_{2})\left( b_{1},b_{2}\right) =\left( a_{1}b_{1}+\gamma \overline{%
b_{2}}a_{2},a_{2}\overline{b_{1}}+b_{2}a_{1}\right) ,
\end{equation*}%
\newline
obtaining an algebra structure over $A\oplus A.$ This algebra, denoted by $%
\left( A,\gamma \right) .$ is called the \textit{algebra obtained from }$A$%
\textit{\ by the Cayley-Dickson process.} $\,$It is clear that $\dim \left(
A,\gamma \right) =2\dim A$.

For $x\in \left( A,\gamma \right) $, $x=\left( a_{1},a_{2}\right) $, the map 
\begin{equation*}
\,\,\,\overline{\phantom{x}}:\left( A,\gamma \right) \rightarrow \left(
A,\gamma \right) \,,\,\,x\rightarrow \bar{x}\,=\left( \overline{a}_{1},\text{%
-}a_{2}\right) ,
\end{equation*}%
\newline
is a scalar involution of the algebra $\left( A,\gamma \right) $, which
extend the involution $\overline{\phantom{x}}\,\,\,$of the algebra $A.$ Let 
\begin{equation*}
\,t\left( x\right) =t(a_{1})
\end{equation*}%
and$\,\,\,$ 
\begin{equation*}
n\left( x\right) =n\left( a_{1}\right) -\gamma n(a_{2})
\end{equation*}%
be $\,\,$the \textit{trace} and the \textit{norm} of the element $x\in $ $%
\left( A,\gamma \right) ,$ respectively.

For $A=K,$ \thinspace applying this process $t$ times, $t\geq 1,\,\,$we
obtain the following algebra over $K,\,\,$%
\begin{equation*}
A_{t}=\left( \frac{\alpha _{1},...,\alpha _{t}}{K}\right) .
\end{equation*}%
The set $\{1,f_{2},...,f_{q}\},q=2^{t}$ generates a basis with the
properties:%
\begin{equation*}
f_{i}^{2}=\alpha _{i}1,\,\,\alpha _{i}\in K,\alpha _{i}\neq 0,\,\,i=2,...,q
\end{equation*}%
and \ 
\begin{equation*}
f_{i}f_{j}=-f_{j}f_{i}=\beta _{ij}f_{k},\,\,\beta _{ij}\in K,\,\,\beta
_{ij}\neq 0,i\neq j,i,j=\,\,2,...q,
\end{equation*}%
$\ \beta _{ij}$ and $f_{k}$ being uniquely determined by $f_{i}$ and $f_{j}.$

If \thinspace 
\begin{equation*}
x\in A_{t},x=x_{1}1+{\sum\limits_{i=2}^{q}}x_{i}f_{i},
\end{equation*}%
the quadratic form $T_{C}:A_{t}\rightarrow K,$ 
\begin{equation*}
T_{C}=<1,\alpha _{1},\alpha _{2},-\alpha _{1}\alpha _{2},\alpha _{3},\text{%
...},\left( -1\right) ^{t+1}({\prod\limits_{i=1}^{t}}\alpha _{i})>=<1,\beta
_{2},\text{...},\beta _{q}>
\end{equation*}%
is called \textit{the trace form }and the quadratic form $T_{P}=T_{C}\mid
_{(A_{t})_{0}}:(A_{t})_{0}\rightarrow K,\ $ 
\begin{equation*}
T_{P}=<\alpha _{1},\alpha _{2},-\alpha _{1}\alpha _{2},\alpha _{3},\text{...}%
,\left( -1\right) ^{t+1}({\prod\limits_{i=1}^{t}}\alpha _{i})>=<\beta _{2},%
\text{...},\beta _{q}>
\end{equation*}%
is called \textit{the pure trace form} of the algebra $A_{t}.$ We remark
that $T_{C}=<1>\perp T_{P}$ (the orthogonal sum of two quadratic forms) and~ 
$n=n_{C}=<1>\perp -T_{P},$ $\ $therefore$\ \ $%
\begin{equation*}
\ n_{C}=<1,-\alpha _{1},-\alpha _{2},\alpha _{1}\alpha _{2},\alpha _{3},%
\text{...},\left( \text{-}1\right) ^{t}({\prod\limits_{i=1}^{t}}\alpha
_{i})>=<1,-\beta _{2},\text{...},-\beta _{q}>.
\end{equation*}

Generally, these algebras $A_{t}$ of dimension $2^{t},\ $obtained by the
Cayley-Dickson process, are not division algebras for all $t\geq 1$.\ There
are some fields on which, if we apply the Cayley-Dickson process, the
resulting algebras $A_{t}\ $are division \ algebras for all $t\geq 1.$ Such
a construction was given by R. B. Brown in [Br; 67], in which he built a
division algebra $A_{t}$ of dimension $2^{t}$ over the power-series field $%
K\{X_{1},X_{2},...,X_{t}\},$ for every $t.$We will shortly recall using
polynomial rings over $K$ and their rational function field instead of
power-series field over $K$ (as it was done by R.B. Brown ).

It is well known that if an algebra $A$ is finite-dimensional, then it is a
division algebra if and only if $A$ does not contain zero divisors (See
[Sc;66]). Starting from this remark, for every $t$ we construct a division
algebra $A_{t}$ over a field $F_{t}.$ We consider $X_{1},X_{2},...,X_{t}$ be 
$t$ algebraically independent indeterminates over the field $K$ and $%
F_{t}=K\left( X_{1},X_{2},...,X_{t}\right) $ be the rational function field$%
. $ For $i=1,...,t,$ we construct the algebra $A_{i}$ over the rational
function field $K\left( X_{1},X_{2},...,X_{i}\right) \,\,\,$by taking $%
\alpha _{j}=X_{j}$ for $j=1,2,...,\,i.\,\,$Let $A_{0}=K.$ $\,$By\thinspace
\thinspace \thinspace \thinspace \thinspace induction over$\,\,\,i,$%
\thinspace \thinspace if we suppose that $A_{i-1}$ is a division algebra
over the field $F_{i-1}=K\left( X_{1},X_{2},...,X_{i-1}\right) $, we may
prove\thinspace \thinspace \thinspace \thinspace that\thinspace \thinspace
\thinspace \thinspace the\thinspace \thinspace \thinspace algebra$%
\,\,\,A_{i} $ is a division algebra over the field $F_{i}=K\left(
X_{1},X_{2},...,X_{i}\right) $.

$\,$Let $A_{F_{i}}^{i-1}=F_{i}\otimes _{F_{i-1}}A_{i-1}.$ For $\alpha
_{i}=X_{i}$ we apply the Cayley-Dickson process to the algebra $%
A_{F_{i}}^{i-1}.$ The resulting algebra, denoted by $A_{i},$ is a division
algebra over the field $F_{i}=K\left( X_{1},X_{2},...,X_{i}\right) \,\,\,$of
dimension $2^{i},~$(see [Fl; 13]).\medskip

\textbf{Cassels-Pfister Theorem.} ([La, Ma;01, p.1823, Theorem 1.3.]) 
\textit{Let }$\varphi ,\psi =<1>\bot \psi ^{^{\prime }}$\textit{\ be two
quadratic forms over a field }$K,\,charK\neq 2$\textit{. If }$\varphi $%
\textit{\ is anisotropic over }$K$\textit{\ and }$\varphi _{K\left( \psi
\right) }$\textit{\ is hyperbolic, then }$\alpha \psi <\varphi $\textit{\
for any scalar represented by }$\varphi .$\textit{\ In particular, }$\dim
\varphi \geq \dim \psi .$\bigskip

\textbf{Springer's Theorem.} ([La, Ma;01, p.1823, Theorem 1.1.]) \textit{Let 
}$\varphi _{1},\ \varphi _{2}$\textit{\ be two quadratic forms over a field }%
$K$\textit{\ and}$~K(X)$\textit{\ be the rational function field \ over }$K.$%
\textit{\ Then, the quadratic form }$\varphi _{1}\perp X\varphi _{2}$\textit{%
\ is isotropic over }$K(X)$\textit{\ if and only if }$\varphi _{1}$\textit{\
or }$\varphi _{2}$\textit{\ is isotropic over }$K.\medskip $

\textbf{Proposition 1.1. }(Theorem 4.1, [Ka, Me; 03])

\textit{Let} $\varphi $ \textit{and} $\phi $ \textit{be two anisotropic
quadratic forms over a field} $K.$ \textit{Assuming that} $\phi $ \textit{is
isotropic over} $K(\varphi ),$ \textit{we get:}

\textit{1)} $dim_{es}(\varphi )\leq $ $dim_{es}(\phi )$.

\textit{2) The equality} $dim_{es}(\varphi )=$ $dim_{es}(\phi )~$\textit{%
holds if and only if} $\varphi $ \textit{is isotropic over} $K(\phi )$%
.\medskip

\textbf{Proposition 1.2.} ([Fl; 11], Proposition 3.5.) \textit{Let }$A$%
\textit{\ be an algebra over a field }$K$\textit{\ obtained by the
Cayley-Dickson process, of dimension }$q=2^{t},T_{C}$\textit{\ ~and }$T_{P}$%
\textit{\ be its trace and pure trace forms}$.\,\,$\textit{Let }$n=2^{k}-1.%
\mathit{\ \,}$\textit{If\thinspace \thinspace \thinspace \thinspace }$t\geq
2\,\,$ \textit{and}\thinspace \thinspace $k>1,$ \textit{then} \underline{$s$}%
$\left( A\right) \leq 2^{k}-1$ \textit{if and only if} $\,\,<1>\bot
(2^{k}-1)\times T_{P}$ \textit{is isotropic}.\medskip

\textbf{Proposition 1.3.} ([Fl; 11], Proposition 3.7.) \textit{Let }$A$%
\textit{\ be an algebra over a field }$K$\textit{\ obtained by the
Cayley-Dickson process, of dimension }$q=2^{t},T_{C}$\textit{\ ~and }$T_{P}$%
\textit{\ be its trace and pure trace forms}$.$ \textit{If \ }$k\geq t,$%
\textit{then }$s\left( A\right) \leq 2^{k}$\textit{\ if and only if the form 
}$\left( 2^{k}+1\right) \times <1>\perp \left( 2^{k}-1\right) \times T_{P}$%
\textit{\ is isotropic}.$\medskip $

\textbf{Proposition 1.4.} ([Kn; 76], Theorem 3.3 and Example 4.1, [La, Ma;
01], Theorem 1.5.)

i) \textit{Let} $K$ \textit{be a field and} $\varphi $ \textit{be a
quadratic form over} $K$ \textit{and} $L$ \textit{be a field extensions of} $%
K.$ \textit{If }$\varphi _{L}$ \textit{is isotropic, then there is a} $K$%
\textit{-place from} $K\left( \varphi \right) $ \textit{to} $L.$

ii) \textit{If }\ $\varphi $ \textit{is a Pfister neighbor of an }$n-$%
\textit{\ fold Pfister form, }$\Psi ,$ \textit{then }$K\left( \varphi
\right) $ \textit{and }$K\left( \Psi \right) ~$\textit{are } $K-$\textit{%
equivalent(that means there exist an }$K-$\textit{place between them)}%
.\medskip

\textbf{Proposition 1.5.} ([Fl; 11], Proposition 3.8) \textit{Let }$K$ 
\textit{be a field and let }$A$\textit{\ be an algebra over the field }$K$%
\textit{\ obtained by the Cayley-Dickson process, of dimension }$q=2^{t}.$%
\textit{\ }

\textit{i) If\thinspace } $k\geq 2,$ \textit{then\ } \underline{$s$}$\left(
A\right) \leq 2^{k}-1$ \textit{if and only if}$\,$ $s\left( A\right) \leq
2^{k}-1.$

\textit{ii) If}\thinspace \thinspace \thinspace \underline{$s$}$\left(
A\right) =n\,$ \textit{and\thinspace } $k\geq 2\,$ \textit{such
that\thinspace } $2^{k-1}\leq n<2^{k},$ \textit{then} $\,s\left( A\right)
\leq 2^{k}-1.$

\textit{iii)} \textit{If} \underline{$s$}$\left( A\right) =1$ \textit{then} $%
\,s\left( A\right) \leq 2.\medskip $

\textbf{Proposition 1.6.} ([Fl; 11], Proposition 3.1)\textit{\ Let }$A$%
\textit{\ be an algebra over a field }$K$\textit{\ obtained by the
Cayley-Dickson process, of dimension }$q=2^{t},T_{C}$\textit{\ and }$T_{P}$%
\textit{\ be its trace and pure trace forms}$.$\textit{The following
statements are true:}

\textit{i) If~~ }$s\left( A\right) \leq n$\textit{\ then }$-1$\textit{\ is
represented by the quadratic form }$n\times T_{C}.$

\textit{ii)} \textit{For }$n\in \mathbb{N}-\{0\},$\textit{\ } \textit{if the
quadratic form }$<1>\perp n\times T_{P}\,\,\,$\textit{is isotropic over }$%
K,\,\,\,$\textit{then }$s(A)\leq n.\medskip $

\begin{equation*}
\end{equation*}%
\textbf{2. Main results}%
\begin{equation*}
\end{equation*}

Let $A_{t}$\ be a division algebra over the field $\ \
K=K_{0}(X_{1},...,X_{t})$\ obtained by the Cayley-Dickson process and
Brown's construction of dimension $q=2^{t},$ with $\ \ K_{0}$\textit{\ } a
formally real field, \textit{\ }$X_{1},...,X_{t}$ algebraically independent
indeterminates over the field \ $K_{0},\ $ $T_{C}$\ and $T_{P}$\ its trace
and pure trace forms$.$ Let%
\begin{equation*}
\varphi _{n}=<1>\bot n\times T_{P},\psi _{m}=<1>\bot m\times T_{C},n\geq 1,
\end{equation*}%
\begin{equation}
A_{t}\left( n\right) =A_{t}\otimes _{K}K\left( <1>\bot n\times T_{P}\right)
,n\in \mathbb{N}-\{0\}.  \tag{2.1.}
\end{equation}

We denote $K_{n}=K\left( <1>\bot n\times T_{P}\right) ,$ and let $%
n_{C}^{A_{t}}$ be the norm form of the algebra $A_{t}$.

\bigskip

\textbf{Remark 2.1. }We have that 
\begin{equation*}
A_{t}\left( n\right) =\left( \frac{X_{1},X_{2},...,X_{t}}{%
K_{0}(X_{1},...,X_{t})}\right)
\end{equation*}%
and, since in [Fl; 13], Proposition 2.3 i), was proved that the norm form $%
n_{C}^{A_{t}\left( n\right) }$ is anisotropic over $K_{n},$ we have that for 
$t\in \{2,3\}$ the \ obtained algebras are division algebras. For $\ t=2,$
the quaternion algebra $\left( \frac{X_{1},X_{2}}{K_{0}(X_{1},X_{2})}\right) 
$ is a division algebra and for $t=3,$ the octonion algebra $\left( \frac{%
X_{1},X_{2},X_{3}}{K_{0}(X_{1},X_{2},X_{3})}\right) $ is a division
algebra.\medskip

\textbf{Proposition 2.2.} \textit{With the above notations,} \textit{for} $%
n\leq 2^{t}-1,$ \textit{the form} $n\times T_{P}$ \textit{is anisotropic over%
} $K_{n}.\medskip $

\textbf{Proof.} \ First of all, we remark that by repeated application of
the Springer's theorem, we have the forms $n\times T_{P}$ \ and $<1>\bot
n\times T_{P}$ anisotropic over $K.$In [Fl; 13], Proposition 2.6, was
obtained that $i_{1}\left( <1>\bot n\times T_{P}\right) =1,$\textit{\ }for
all\textit{\ \ }$n\in \mathbb{N}-\{0\},$\textit{\ }where\textit{\ }$T_{P}$%
\textit{\ }is the pure trace form for the algebra\textit{\ }$A_{t},t\geq 2.$
Now, using Proposition 1.1, we have $dim_{es}\left( \varphi _{n}\right)
=~dim\left( \varphi _{n}\right) .$ It results that $dim\left( \varphi
_{n}\right) =n\left( 2^{t}-1\right) +1$ and $dim(n\times T_{P})=n\left(
2^{t}-1\right) .$ Since $dim_{es}\left( n\times T_{P}\right) \leq
~dim(n\times T_{P})\,<~dim\left( \varphi _{n}\right) =~dim_{es}\left(
\varphi _{n}\right) ,$ we have the form $n\times T_{P}$ anisotropic over $%
K_{n}.\Box \smallskip $

\textbf{Proposition 2.3.} \textit{\ For a fixed }$t$ \textit{and for each} $%
n\leq 2^{t}-1,$ \textit{it results that }$s\left( A_{t}\left( n\right)
\right) =\underline{s}\left( A_{t}\left( n\right) \right) =n.\medskip $

\textbf{Proof.} \ From [Fl; 13], Theorem 2.9, we know that $s\left(
A_{t}\left( n\right) \right) =n$ and $\underline{s}\left( A_{t}\left(
n\right) \right) \in \{n-1,n\},$ for $n\leq 2^{t}-1.$ Let $q=2^{t}$ and $%
\{1,f_{2},...,f_{q}\}$ be a basis in $A_{t}\left( n\right) $. If we consider 
$\underline{s}\left( A_{t}\left( n\right) \right) =n-1,$ it results that we
can find the nonzero elements 
\begin{equation*}
u_{i}={x_{i1}+}x_{i2}f_{2}+...+x_{iq}f_{q}\in A_{t}\left( n\right) ,~
\end{equation*}%
with%
\begin{equation*}
u_{i}^{\prime \prime }=x_{i2}f_{2}+...+x_{iq}f_{q}\in A_{t}\left( n\right) ,
\end{equation*}%
the pure part of $u_{i},$ where $x_{ij}\in K,i\in \{1,2,...,n\},j\in
\{1,2,...,q\},q=2^{t},$ such that $0=u_{1}^{2}+...+u_{n}^{2}.$ We know that
the pure elements in $A_{t}\left( n\right) $ form a vector subspace over $%
K_{n}$, denoted with $A_{t}\left( n\right) ^{\ast }.$ We obtain\ 
\begin{equation*}
\sum\limits_{i=1}^{n}(x_{i1}^{2}+(u_{i}^{\prime \prime
})^{2}+2x_{i1}u_{i}^{\prime \prime })=0,
\end{equation*}%
therefore \ 
\begin{equation}
\sum\limits_{i=1}^{n}x_{i1}^{2}+\sum\limits_{i=1}^{n}(u_{i}^{\prime \prime
})^{2}=0  \tag{2.2.}
\end{equation}%
\ and%
\begin{equation}
\ \sum\limits_{i=1}^{n}x_{i1}u_{i}^{\prime \prime }=0.  \tag{2.3.}
\end{equation}

\textbf{Case 1.} If $x_{i1}=0,$ $\forall $ $i\in \{1,2,...,n\}.$ It results
that 
\begin{equation*}
\sum\limits_{i=1}^{n}(u_{i}^{\prime \prime })^{2}=0,
\end{equation*}%
hence, it follows that \ $n\times T_{P}$ is isotropic over $K_{n}$, false,
using Proposition 2.2.

\textbf{Case 2. }If there is an element $x_{i1}\neq 0,$ from relation $%
\left( 2.3\right) ,$ we obtain that \textbf{\ }there is a $K_{n}$ vector
subspace of \ $A_{t}\left( n\right) ^{\ast }$ of dimension $n-1,$ denoted $%
V_{n-1},$ which contains $u_{i}^{\prime \prime },i\in \{1,2,...,n\}.$ We
consider 
\begin{equation*}
\eta :V_{n-1}\rightarrow K_{n},\eta \left( u\right) =u^{2}.
\end{equation*}

It results that $\eta $ is a subform of $T_{p}$ of \ dimension \ at most $%
n-1 $ and, from relation $\left( 2.2\right) ,$ we have that $n\times
(<1>\bot \eta )$ is isotropic over $K_{n}.$ We have that the form $<1>\bot
n\times T_{P}$ are anisotropic over $K$ and $i_{1}\left( <1>\bot n\times
T_{P}\right) =1.$ Computing dimension of the form $n\times (<1>\bot \eta ),$
we obtain the following result: 
\begin{equation*}
dim\left( n\times (<1>\bot \eta )\right) \leq n^{2}\,<n\left( 2^{t}-1\right)
+1=dim\left( <1>\bot n\times T_{P}\right) .
\end{equation*}%
Therefore, using again Proposition 2.2, we have that $dim_{es}\left( \varphi
_{n}\right) =dim\left( \varphi _{n}\right) .$ Therefore $dim_{es}\left(
n\times (<1>\bot \eta )\right) \leq dim(n\times (<1>\bot \eta ))\,<dim\left(
\varphi _{n}\right) =dim_{es}\left( \varphi _{n}\right) ,$ then the form $%
~n\times (<1>\bot \eta )$ is anisotropic over $K_{n},$ which it is a
contradiction. It results that $\underline{s}\left( A_{t}\left( n\right)
\right) =n.\Box \smallskip $

The above proposition generalized to algebras obtained by the Cayley-Dickson
process Theorem 3.3 and Theorem 3.4 from [O' Sh; 06].\medskip

\textbf{Theorem 2.4.} \textit{With the above notation, for each }$n\in 
\mathbb{N}-\{0\}$\textit{\ there is an algebra }$A_{t}\left( n\right) $%
\textit{\ such that }$s\left( A_{t}\left( n\right) \right) =$ \underline{$s$}%
$\left( A_{t}\left( n\right) \right) =n.\medskip $

\textbf{Proof.} Let $n\in \mathbb{N}-\{0\}$ \ and \ $t\ \ $be$\ $the least
positive integer such that $n\leq 2^{t}.$ For $n=2^{t},$ the result was
given in [Fl; 13], Example 2.5. We assume that $n<2^{t}$ and we apply
Proposition 2.3$.\Box \smallskip $

\textbf{Theorem 2.5.} \textit{1)} \textit{For} $n=2^{k},~k\geq 2,$ \textit{a
natural number, and } $t\geq 3,~$\textit{we have that} $s\left( A_{t}\left(
n\right) \right) =$ \underline{$s$}$\left( A_{t}\left( n\right) \right)
=s\left( A_{t-1}\left( n\right) \right) =$ \underline{$s$}$\left(
A_{t-1}\left( n\right) \right) =n.$

\textit{2) For} $n=2^{k}+1,$ $k$ \textit{a natural number, and } $t\geq 3,~$%
\textit{we have that} $s\left( A_{t}\left( n\right) \right) =s\left(
A_{t-1}\left( n\right) \right) =n.$

\textit{3) For each natural number} $\,n,0<n<2^{t-1},$ $t\geq 3,$ \textit{we
have that} $s\left( A_{t}\left( n\right) \right) =s\left( A_{t-1}\left(
n\right) \right) =n.\medskip $

\textbf{Proof.} \ 1) From the above, we have $<1>\bot n\times T_{P}$ is
anisotropic over $K,$ for all $n.$ If we suppose that \underline{$s$}$\left(
A_{t-1}\left( n\right) \right) \leq 2^{k}-1,$ from Proposition 1.2, we have
that $<1>\bot \left( 2^{k}-1\right) \times T_{P}$ is isotropic over $%
K_{2^{k}},$ with $\varphi _{n}=\varphi _{2^{k}}=<1>\bot 2^{k}\times T_{P}$.
From [Fl; 13], Proposition 2.6, we have that $i_{1}\left( <1>\bot n\times
T_{P}\right) =1,$ for all $n\in \mathbb{N}-\{0\}.$ It results the following
relation $dim_{es}\left( <1>\bot n\times T_{P}\right) =dim\left( <1>\bot
n\times T_{P}\right) .$ Since 
\begin{equation*}
dim\left( <1>\bot 2^{k}\times T_{P}\right) >dim\left( <1>\bot \left(
2^{k}-1\right) \times T_{P}\right) ,
\end{equation*}%
from Proposition 1.1, we have that $<1>\bot \left( 2^{k}-1\right) \times
T_{P}$ is anisotropic over $K_{2^{k}},$a contradiction. Therefore, using
Example 2.5 from [Fl; 13], we obtain that \underline{$s$}$\left(
A_{t-1}\left( n\right) \right) =s\left( A_{t-1}\left( n\right) \right) =n.$

2) Using Proposition 1.3, if $s\left( A_{t}\left( n\right) \right) \leq
2^{k},$ then $\left( 2^{k}+1\right) \times <1>\bot \left( 2^{k}-1\right)
\times T_{P}$ is isotropic over $K_{2^{k}}.$ We have that the form $%
X_{t}\left( 2^{k}+1\right) \times n_{C}^{A_{t-1}\left( n\right) }$is a
subform of $<1>\bot \left( 2^{k}+1\right) \times T_{P}.$ From here, it
results that $\left( 2^{k}+1\right) \times <1>\bot \left( 2^{k}-1\right)
\times T_{P}$ is isotropic over $K\left( X_{t}\left( 2^{k}+1\right) \times
n_{C}^{A_{t-1}\left( n\right) }\right) $. Since 
\begin{equation*}
K\left( X_{t}\left( 2^{k}+1\right) \times n_{C}^{A_{t-1}\left( n\right)
}\right) \simeq K\left( \left( 2^{k}+1\right) \times n_{C}^{A_{t-1}\left(
n\right) }\right) ,
\end{equation*}%
from [Sch; 85], Remarks 5.2, ii), p.154, we have that $\left( 2^{k}+1\right)
\times <1>\bot \left( 2^{k}-1\right) \times T_{P}$ is isotropic over $%
K\left( \left( 2^{k}+1\right) \times n_{C}^{A_{t-1}\left( n\right) }\right)
. $ By repeated used of of Springer's Theorem, we have that $\left(
2^{k}+1\right) \times <1>\bot \left( 2^{k}-1\right) \times <X_{1}>$ or $%
\left( 2^{k}-1\right) \times <1,-X_{1}>$ is isotropic over $K\left( \left(
2^{k}+1\right) \times n_{C}^{A_{t-1}\left( n\right) }\right) .$ Since 
\begin{equation*}
dim\left( \left( 2^{k}+1\right) \times <1>\bot \left( 2^{k}-1\right) \times
<X_{1}>\right) \text{=}2^{k+1}\leq dim\left( \left( 2^{k}+1\right) \times
n_{C}^{A_{t-1}\left( n\right) }\right)
\end{equation*}%
and 
\begin{equation*}
dim\left( \left( 2^{k}-1\right) \times <1,-X_{1}>\right) \text{=}%
2^{k+1}-2<2^{k+1}\leq dim\left( \left( 2^{k}+1\right) \times
n_{C}^{A_{t-1}\left( n\right) }\right) ,
\end{equation*}%
using Theorem 1 from [Hoff; 95], we have that $\left( 2^{k}+1\right) \times
<1>\bot \left( 2^{k}-1\right) \times <X_{1}>$ and $\left( 2^{k}-1\right)
\times <1,-X_{1}>$ remain anisotropic over $K\left( \left( 2^{k}+1\right)
\times n_{C}^{A_{t-1}\left( n\right) }\right) ,$ which is false. Therefore $%
s\left( A_{t}\left( n\right) \right) =2^{k}+1.$

3) Since $[\frac{n}{2^{t}}]=[\frac{n}{2^{t-1}}]=0,$ where $[]$ denote the
integer part, using Theorem 2.9, from [Fl; 13], we obtain the asked result.$%
\Box \smallskip $

The above Theorem, i) and ii) generalized to algebras obtained by the
Cayley-Dickson process Theorems 3.2 and 3.4 from [O' Sh; 10].\medskip

\textbf{Theorem 2.6.} \textit{Let} $n\in \mathbb{N},n=2^{k}r,r$ \ \textit{be}
\textit{an odd number. If the form} $\left( 2^{k}+1\right) \times <1>\bot
\left( n-1\right) \times T_{P}$ \textit{is isotropic over} $K$\textit{,} 
\textit{then} $s\left( A_{t}\right) \leq n.\medskip $

\textbf{Proof.} If the form $\left( 2^{k}+1\right) \times <1>\bot \left(
n-1\right) \times T_{P}$ is isotropic over $K,$ therefore it is universal.
Writing the form $\left( 2^{k}+1\right) \times <1>\bot \left( n-1\right)
\times T_{P}=2^{k}\times <1>\bot \left( <1>\bot \left( n-1\right) \times
T_{P}\right) ,$ we obtain that there are the elements $U\in A_{t},V\in K,$
such that $2^{k}\times <1>\left( U\right) =V$ and $(<1>\bot \left(
n-1\right) \times T_{P})\left( U\right) =-V.$ It results that%
\begin{equation*}
-V=\alpha
^{2}+X_{1}\sum\limits_{i=1}^{n-1}x_{i1}^{2}+X_{2}\sum%
\limits_{i=1}^{n-1}x_{i2}^{2}-X_{1}X_{2}\sum\limits_{i=1}^{n-1}x_{i3}^{2}+
\end{equation*}%
\begin{equation*}
+...\left( \text{-}1\right) ^{t+1}({\prod\limits_{i=1}^{t}}%
X_{i})\sum\limits_{i=1}^{n-1}x_{iq}^{2},q=2^{t},x_{ij}\in K.
\end{equation*}

\bigskip We obtain%
\begin{equation*}
-1\text{=}\frac{1}{V^{2}}(\alpha ^{2}V\text{+}X_{1}V\sum%
\limits_{i=1}^{n-1}x_{i1}^{2}\text{+}X_{2}V\sum\limits_{i=1}^{n-1}x_{i2}^{2}-
\end{equation*}

\begin{equation*}
\text{--}X_{1}X_{2}V\sum\limits_{i=1}^{n-1}x_{i3}^{2}\text{+...}\left( \text{%
-}1\right) ^{t+1}({\prod\limits_{i=1}^{t}}X_{i})V\sum%
\limits_{i=1}^{n-1}x_{iq}^{2}).
\end{equation*}
\newline

\bigskip Since $(2^{k}\times <1>)\left( U\right) =V,$ we have that 
\begin{equation*}
V=\sum\limits_{i=1}^{2^{k}}y_{1i}^{2},y_{1i}\in K.
\end{equation*}

We will prove that there are \ elements $y_{2i}\in K,$ such that 
\begin{equation*}
\sum\limits_{i=1}^{n}y_{2i}^{2}=V\sum\limits_{i=1}^{n-1}x_{i1}^{2}
\end{equation*}%
and

\begin{equation*}
\sum\limits_{i=1}^{n}y_{1i}y_{2i}=0,
\end{equation*}%
completing with zero if $2^{k}<n.$ If all \ $y_{1i}=0,$ then we put all $%
y_{2i}=0.$ If there are $y_{1i}\neq 0,$ we have that $<V,V\sum%
\limits_{i=1}^{n-1}x_{i1}^{2}>\simeq V<1,\sum\limits_{i=1}^{n-1}x_{i1}^{2}>$
is a subform of \ $Vn\times <1>\simeq n\times <1>,$ since $%
\sum\limits_{i=1}^{n-1}x_{i1}^{2}$ is represented by the form $\left(
n-1\right) \times <1>$ and $V2^{k}\times <1>\simeq 2^{k}\times <1>.$
Therefore such elements exist. In the same way, we can find $y_{3i}\in K,$
such that 
\begin{equation*}
\sum\limits_{i=1}^{n}y_{3i}^{2}=V\sum\limits_{i=1}^{n-1}x_{i2}^{2}
\end{equation*}%
and

\begin{equation*}
\sum\limits_{i=1}^{n}y_{1i}y_{3i}=0,etc.
\end{equation*}

\bigskip It results that 
\begin{equation*}
\sum\limits_{i=1}^{n}\left( \frac{\alpha y_{1i}}{V}+\frac{y_{2i}}{V}f_{2}+%
\frac{y_{3i}}{V}f_{3}+...+\frac{y_{qi}}{V}f_{q}\right) ^{2}=
\end{equation*}

\begin{equation*}
=\frac{1}{V^{2}}(\alpha ^{2}V\text{+}X_{1}V\sum\limits_{i=1}^{n-1}x_{i1}^{2}%
\text{+}X_{2}V\sum\limits_{i=1}^{n-1}x_{i2}^{2}-
\end{equation*}

\begin{equation*}
-X_{1}X_{2}V\sum\limits_{i=1}^{n-1}x_{i3}^{2}\text{+...}\left( \text{-}%
1\right) ^{t+1}({\prod\limits_{i=1}^{t}}X_{i})V\sum%
\limits_{i=1}^{n-1}x_{iq}^{2})\text{=}-1,
\end{equation*}%
where $\{1,f_{2},...,f_{q}\}$ is a basis in $A_{t}.$ From here, we have $%
s\left( A_{t}\right) \leq n.\Box \smallskip $

The above result generalized to algebras obtained by the Cayley-Dickson
process Theorem 3.11 from [O' Sh; 10].\medskip

\textbf{Definition 2.7.} Let $n\in \mathbb{N},n=2^{k}r,r$ \ be an odd
number, $n>1,$ $\varphi _{n}^{\prime }=\left( 2^{k}+1\right) \times <1>\bot
\left( n-1\right) \times T_{P},$ be the above form and 
\begin{equation}
A_{t}^{\prime }\left( n\right) =A_{t}\otimes _{K}K\left( \left(
2^{k}+1\right) \times <1>\bot \left( n-1\right) \times T_{P}\right) ,n\in 
\mathbb{N}-\{0\}.  \tag{2.4.}
\end{equation}

We denote $K_{n}^{\prime }=K\left( \left( 2^{k}+1\right) \times <1>\bot
\left( n-1\right) \times T_{P}\right) ,$ and let $n_{C}^{A_{t}^{\prime }}$
be the norm form of the algebra $A_{t}^{\prime }$.\medskip

\textbf{Proposition 2.8.} \textit{The norm form} $n_{C}^{A_{t}^{\prime
}\left( n\right) }$ \textit{is anisotropic over} $K_{n}^{\prime }.\medskip $

\textbf{Proof.} \ For $n>1$ and we suppose that $n_{C}^{A_{t}^{\prime
}\left( n\right) }$ is isotropic over $K_{n}^{\prime }$. We have that $%
n_{C}^{A_{t}^{\prime }\left( n\right) }$ is a Pfister form. A Pfister form
is isotropic if and only if it is hyperbolic, therefore $n_{C}^{A_{t}^{%
\prime }\left( n\right) }\,$\ is hyperbolic over $K_{n}^{\prime }$. From
Brown's construction, described above, the algebra $A_{t}$ is a division
algebra, therefore $n_{C}^{A_{t}}$ is anisotropic. Since an isotropic form
is universal, using Cassels-Pfister Theorem, we have that for each element $%
a\in K_{n}^{\ast },$ the form $a$ $\varphi _{n}^{\prime }$ is a subform of
the norm form $n_{C}^{A_{t}^{\prime }\left( n\right) }.$ We have $%
dim~\varphi _{n}^{\prime }=2^{k}+1+\left( n-1\right) \left( 2^{t}-1\right) $
and $dim~n_{C}^{A_{t}^{\prime }\left( n\right) }=2^{t}$ and, from here, $%
dim~\varphi _{n}^{\prime }>dim~n_{C}^{A_{t}^{\prime }\left( n\right) },$
which is false. $\Box \smallskip $ \medskip

\textbf{Remark 2.9. }If the algebra $A$ is an algebra obtained by the
Cayley-Dickson process, of dimension greater than $2$ and if $\,n_{C\text{ }%
} $ is isotropic, then $s\left( A\right) =\underline{s}\left( A\right) =1.$
If $-1\notin K^{\ast 2},$ since $n_{C}=<1>\perp -T_{P}$ and $n_{C\text{ }}$
is a Pfister form, we obtain that $-T_{P}$ is isotropic, therefore $T_{C}$
is isotropic and, from Proposition 1.2, we have that $s\left( A\right) =%
\underline{s}\left( A\right) =1.$Therefore if the form $n_{C}^{A_{t}^{\prime
}\left( n\right) }$ $\ $is anisotropic, then $A_{t}^{\prime }\left( n\right) 
$ has level and sublevel greater than $1.\medskip $

\textbf{Proposition 2.10.} \textit{We consider} $k\geq 1,$ \textit{a natural
number. Then the quadratic forms} $2^{k}\times <1>\bot \left( 2^{k}-1\right)
\times T_{P}^{t-1}$ \textit{and} $\left( 2^{k}-1\right) \times
n_{C}^{A_{t-1}^{\prime }\left( n\right) }$ \textit{are anisotropic over} $%
K\left( \left( 2^{k}+1\right) \times <1>\right) $.\medskip

\textbf{Proof.} Supposing that the forms $2^{k}\times <1>\bot \left(
2^{k}-1\right) \times T_{P}^{t-1}$ and $\left( 2^{k}-1\right) \times
n_{C}^{A_{t-1}^{\prime }\left( n\right) }$ are isotropic over $K\left(
\left( 2^{k}+1\right) \times <1>\right) ,$ from Springer's Theorem, we have
that the quadratic forms $2^{k}\times <1>$ and $\left( 2^{k}-1\right) \times
<1>$ are isotropic. It is clear that these forms are Pfister neighbors of
the Pfister form $2^{k}\times <1>$. From here, using Proposition 1.4, ii),
we have that their functions fields are $K-$equivalent to $K(2^{k}\times
<1>).$ If the forms $2^{k}\times <1>$ or $\left( 2^{k}-1\right) \times <1>$
are isotropic over $K\left( \left( 2^{k}+1\right) \times <1>\right) $, using
Proposition 1.4, i), we have that there is a $K-$place from $K(2^{k}\times
<1>)$ to $K\left( \left( 2^{k}+1\right) \times <1>\right) $. From here,
since $2^{k}\times <1>$ is isotropic over $K\left( 2^{k}\times <1>\right) ,$
it results that it is isotropic over $K\left( \left( 2^{k}+1\right) \times
<1>\right) ,~$therefore it is hyperbolic over $K\left( \left( 2^{k}+1\right)
\times <1>\right) .$ Since 
\begin{equation*}
dim(2^{k}\times <1>)=2^{k}<dim((2^{k}+1)\times <1>)=2^{k}+1,
\end{equation*}%
we obtain a contradiction with Cassels-Pfister Theorem. $\Box \smallskip $

The above proposition generalized Proposition 3.3 from [La, Ma; 01] to
algebras given in Definition 2.7.\medskip\ \medskip

\textbf{Proposition 2.11.} \textit{With the above notations, the form} $%
2^{k}\times <1>\bot (2^{k}-1)\times T_{P}^{t}$ \textit{is anisotropic over} $%
K_{2^{k}}^{\prime }.\medskip $

\textbf{Proof.} \ Since $(2^{k}+1)\times <1>$ is a subform of the form%
\newline
$\left( 2^{k}+1\right) \times <1>\bot (2^{k}-1)\times T_{P}^{t},$ we have
that there is a \ $K-$place from $K\left( \left( 2^{k}+1\right) \times
<1>\bot (2^{k}-1)\times T_{P}^{t}\right) $ to $K\left( (2^{k}+1)\times
<1>\right) .$ If $2^{k}\times <1>\bot (2^{k}-1)\times T_{P}^{t}$ is
isotropic over \newline
$K\left( \left( 2^{k}+1\right) \times <1>\bot (2^{k}-1)\times
T_{P}^{t}\right) ,$ then we have \ that $2^{k}\times <1>\bot (2^{k}-1)\times
T_{P}^{t}$ is isotropic over $K\left( (2^{k}+1)\times <1>\right) .$ From
Springer's Theorem, we have that the form $2^{k}\times <1>\bot \left(
2^{k}-1\right) \times T_{P}^{t-1}$ or $\left( 2^{k}-1\right) \times
n_{C}^{A_{t-1}^{\prime }\left( n\right) }$ is isotropic over $K\left( \left(
2^{k}+1\right) \times <1>\right) ,$ false if we use Proposition 2.10.$\Box
\smallskip $ \medskip

The above proposition generalized Proposition 3.3 from [La, Ma; 01] and
Proposition 3.5 from [Pu; 05] to algebras obtained by the Cayley-Dickson
process given in Definition 2.7.\medskip\ \medskip

\textbf{Theorem 2.12.} \textit{For} $n=2^{k},~k\geq 2,$ \textit{a natural
number, algebras} $A_{t}^{\prime }\left( n\right) ,$ \textit{given in
relation} $\left( 2.4\right) ,$ \textit{have the level and sublevel equal
with} $n$.$\medskip $

\textbf{Proof.} For the level case, from Theorem 2.6, we have that $%
s(A_{t}^{\prime }\left( n\right) )\leq 2^{k}.$ If we have $s(A_{t}^{\prime
}\left( n\right) )<2^{k},$ then the quadratic form $2^{k}\times <1>\bot
(2^{k}-1)\times T_{P}^{t}$ is isotropic over $K_{2^{k}}^{\prime },$ which is
a contradiction with the Proposition 2.11.

For the sublevel case, using Proposition 1.5, if \underline{$s$}$\left(
A_{t}^{\prime }\left( n\right) \right) \leq 2^{k}-1$\ we have that $%
s(A_{t}^{\prime }\left( n\right) )\leq 2^{k}-1,$ false.$\Box \smallskip $

The above proposition generalized Theorem 3.4 from [Pu; 05] to algebras
obtained by the Cayley-Dickson process given in Definition 2.7.\medskip

\textbf{Theorem 2.13.} \ \textit{We have that} $i_{1}\left( \left(
2^{k}+1\right) \times <1>\bot \left( n-1\right) \times T_{P}^{t}\right) =1,$%
\textit{\ for all \ }$n\in \mathbb{N}-\{0\},$\textit{\ where }$T_{P}$\textit{%
\ is the pure trace form for the algebra }$A_{t}^{\prime },t\geq 2.\medskip $

\textbf{Proof.} \ We use induction after $t.$ For $t=2,$ the result was
proved in [O' Sh; 10], Theorem 3.13. We suppose that the result is true for $%
t-1$ and we will prove for $t.$ Assuming that $i_{1}\left( \left(
2^{k}+1\right) \times <1>\bot \left( n-1\right) \times T_{P}^{t}\right) >1,$
we obtain that $2^{k}\times <1>\bot \left( n-1\right) \times T_{P}$ is
isotropic over \ $K\left( \left( 2^{k}+1\right) \times <1>\bot \left(
n-1\right) \times T_{P}\right) ,$ since%
\begin{equation*}
dim\left( \left( 2^{k}+1\right) \times <1>\bot \left( n-1\right) \times
T_{P}^{t}\right) -dim\left( 2^{k}\times <1>\bot \left( n-1\right) \times
T_{P}^{t}\right) =1.
\end{equation*}%
Let $\alpha =\left( 2^{k}+1\right) \times <1>\bot \left( n-1\right) \times
T_{P}^{t-1}$ be a subform of the form $\left( 2^{k}+1\right) \times <1>\bot
\left( n-1\right) \times T_{P}^{t},$ therefore $\left( 2^{k}+1\right) \times
<1>\bot \left( n-1\right) \times T_{P}^{t}$ is isotropic over $K\left(
\left( 2^{k}+1\right) \times <1>\bot \left( n-1\right) \times
T_{P}^{t-1}\right) .$ From Proposition 1.4, we have that there is an $K-$%
place from $K\left( \left( 2^{k}+1\right) \times <1>\bot \left( n-1\right)
\times T_{P}^{t}\right) $ to \newline
$K\left( \left( 2^{k}+1\right) \times <1>\bot \left( n-1\right) \times
T_{P}^{t-1}\right) .$ It results that $2^{k}\times <1>\bot \left( n-1\right)
\times T_{P}^{t}$ is isotropic over \newline
$K\left( \left( 2^{k}+1\right) \times <1>\bot \left( n-1\right) \times
T_{P}^{t-1}\right) .$ From Springer's Theorem, we have that $2^{k}\times
<1>\bot \left( n-1\right) \times T_{P}^{t-1}$ or $\left( n-1\right) \times
n_{C}^{A_{t-1}^{\prime }\left( n\right) }$ is isotropic over\newline
$K\left( \left( 2^{k}+1\right) \times <1>\bot \left( n-1\right) \times
T_{P}^{t-1}\right) .$

\textbf{Case 1. }If \ $\left( n-1\right) \times n_{C}^{A_{t-1}^{\prime
}\left( n\right) }$ is isotropic over \newline
$K\left( \left( 2^{k}+1\right) \times <1>\bot \left( n-1\right) \times
T_{P}^{t-1}\right) ,$ therefore for each $r$ such that $2^{r}>n-1,$ we have
that $2^{r}\times n_{C}^{A_{t-1}^{\prime }\left( n\right) }$ is hyperbolic
over $K\left( \left( 2^{k}+1\right) \times <1>\bot \left( n-1\right) \times
T_{P}^{t-1}\right) .$ From Cassels--Pfister Theorem, we have that the form \ 
$\left( 2^{k}+1\right) \times <1>\bot \left( n-1\right) \times T_{P}^{t-1}$
is a subform of the form $2^{r}\times n_{C}^{A_{t-1}^{\prime }\left(
n\right) },$ false.

\textbf{Case 2.} Therefore $2^{k}\times <1>\bot \left( n-1\right) \times
T_{P}^{t-1}$ is isotropic over\newline
$K\left( \left( 2^{k}+1\right) \times <1>\bot \left( n-1\right) \times
T_{P}^{t-1}\right) .$

We consider the form$~~2^{k}\times <1>\bot (n-1)\times T_{P}^{t-1}$ and $P$
be an arbitrary ordering over $K$ such that $X_{1},...X_{n-1}<_{P}0.$ We
remark that such an ordering always exists. Indeed, since $2^{k}\times
<1>\bot (n-1)\times T_{P}^{t-1}$ is anisotropic over $K$, it follows that $%
P_{0}=\{a~\mid ~a=0$ or $a$ is represented by $2^{k}\times <1>\bot
(n-1)\times T_{P}^{t-1}$ $\}$ is a $q-$preordering, therefore there is a $q-$%
ordering $P$ containing $P_{0\text{ }}$ or $-P_{0}.$ \ From here, it results
that $i_{1}\left( 2^{k}\times <1>\bot (n-1)\times T_{P}^{t-1}\right) =2^{k}.$

From induction hypothesis, we have that $i_{1}\left( \left( 2^{k}+1\right)
\times <1>\bot \left( n-1\right) \times T_{P}^{t-1}\right) =1.$ It results
that 
\begin{equation}
dim_{es}\left( \left( 2^{k}\text{+}1\right) \times <1>\bot \left( n\text{-}%
1\right) \times T_{P}^{t-1}\right) >dim_{es}\left( 2^{k}\times <1>\bot
\left( n\text{-}1\right) \times T_{P}^{t-1}\right) .  \tag{2.5.}
\end{equation}

Using Proposition 1.1, we have that $2^{k}\times <1>\bot \left( n-1\right)
\times T_{P}^{t-1}$ is anisotropic over $K\left( \left( 2^{k}+1\right)
\times <1>\bot \left( n-1\right) \times T_{P}^{t-1}\right) ,$ which is a
contradiction. $\Box \medskip $

\textbf{Remark 2.14.} i) If $A$ is an algebra obtained by the Cayley-Dickson
process of dimension \thinspace $q=2^{t},$ such that ~~ \underline{$s$}$%
\left( A\right) <m,$\ \ then the quadratic form $m\times T_{C}$ is
isotropic. Indeed, using Proposition 1.6. i), \ if the quadratic form \ $%
(m+1)\times T_{C}$ is anisotropic, then \underline{$s$}$\left( A\right) \geq
m+1.$

ii) If \ $<1>\perp n\times T_{P}\,\,\,$is isotropic over $K,\,\,\,$then%
\textit{\ }$s(A)\leq n,$ from Proposition 1.6. i).\medskip

\textbf{Remark 2.15.} Let $m=2^{r}s,$ with $s$ an odd number and $n\geq
m,n=2^{k}s^{\prime },s^{\prime }$ an odd number. We remark that for $%
m=2^{r}s,$ we can write $m\times T_{C}\simeq 2^{r}\times <1>\otimes \left(
s\times T_{C}\right) .$ With the above notations, choosing $P$ an arbitrary
ordering over $K$ such that $X_{1},...X_{n}<_{P}0,$ and using Theorem 1.4
from [El, La; 72], we obtain $i_{1}\left( m\times T_{C}\right) \geq 2^{r}.$ $%
\ $\ We have that $dim_{es}\left( \left( 2^{k}+1\right) \times <1>\bot
\left( n-1\right) \times T_{P}^{t}\right) =dim\left( \left( 2^{k}+1\right)
\times <1>\bot \left( n-1\right) \times T_{P}^{t}\right) $ and $%
dim_{es}\left( m\times T_{C}\right) \leq dim\left( m\times T_{C}\right)
-2^{r}+1.$ If \ 
\begin{equation}
dim\left( \left( 2^{k}+1\right) \times <1>\bot \left( n-1\right) \times
T_{P}^{t}\right) >dim\left( m\times T_{C}\right) -2^{r}+1,  \tag{2.6.}
\end{equation}%
\newline
we have that $m\times T_{C}$ is anisotropic over $K\left( \left(
2^{k}+1\right) \times <1>\bot \left( n-1\right) \times T_{P}^{t}\right) ,$
therefore \underline{$s$}$\left( A_{t}^{\prime }\left( n\right) \right) \geq
l,$where $l$ is the least value of $m$ for which relation $\left( 2.6\right) 
$ holds. Since \underline{$s$}$\left( A_{t}^{\prime }\left( n\right) \right)
\leq s\left( A_{t}^{\prime }\left( n\right) \right) ,$ we have that $s\left(
A_{t}^{\prime }\left( n\right) \right) \geq l.$

Using Remark 2.14, i), if $<1>\perp n\times T_{P}~$is isotropic over\newline
$K\left( \left( 2^{k}+1\right) \times <1>\bot \left( n-1\right) \times
T_{P}^{t}\right) ,$ we have that $s\left( A_{t}^{\prime }\left( n\right)
\right) \leq n,$ therefore \underline{$s$}$\left( A_{t}^{\prime }\left(
n\right) \right) \leq n.$ We obtain that 
\begin{equation*}
\underline{s}\left( A_{t}^{\prime }\left( n\right) \right) ,s\left(
A_{t}^{\prime }\left( n\right) \right) \in \lbrack l,n].
\end{equation*}

\textbf{Theorem 2.16.} \textit{With the above notations, for the algebra} $%
A_{t}^{\prime }\left( n\right) $ \textit{of dimension }$2^{t},$ \textit{%
taking} $n=2^{t+2}+2=2\left( 2^{t+1}+1\right) ,m=2^{t+2},$\textit{we have
that}%
\begin{equation*}
\mathit{~}\underline{s}\left( A_{t}^{\prime }\left( n\right) \right)
,s\left( A_{t}^{\prime }\left( n\right) \right) \in \lbrack
2^{t+2},2^{t+2}+2].\medskip
\end{equation*}

\textbf{Proof.} For $n=2^{t+2}+2=2\left( 2^{t+1}+1\right) ,~$we compute%
\newline
$dim\left( \left( 2^{k}+1\right) \times <1>\bot \left( n-1\right) \times
T_{P}^{t}\right) .$ In this situation, $k=1.$ We get\newline
$dim\left( \left( 2^{k}+1\right) \times <1>\bot \left( n-1\right) \times
T_{P}^{t}\right) =3+\left( 2^{t+2}+1\right) \left( 2^{t}-1\right) =$\newline
$=2^{2t+2}-2^{t+2}+2^{t}+2.$

For $m=2^{t+2},$ we compute $dim\left( m\times T_{C}\right) -2^{r}+1.\,\ $In
this situation, $r=t+2.$ We get $dim\left( m\times T_{C}\right)
-2^{r}+1=2^{t+2}2^{t}-2^{t+2}+1,$ therefore relation $\left( 2.6\right) $ is
true.$\Box \medskip $

\textbf{Remark 2.17.} The above Theorem improves the level and sublevel
bounds for the algebras obtained by the Cayley-Dickson process, when the
prescribed level and sublevel are greater than the dimension of the algebra.
Indeed, for $n=34,$ using Theorem 2.7 and Theorem 2.8, from \ [Fl; 13], for
the octonions, we obtain that $s(A_{3}\left( 34\right) )\in \lbrack 30,34]$
and \underline{$s$}$(A_{3}\left( 34\right) )\in \lbrack 29,34].$ From the
above result, we have $s(A_{3}^{\prime }\left( 34\right) ),\underline{s}%
(A_{3}^{\prime }\left( 34\right) )\in \lbrack 32,34],$ therefore the
octonion algebra given by the relation $\left( 2.4\right) $ has better
bounds than octonion algebra given by the relation $\left( 2.1\right) $ for
its level and its sublevel$.\medskip $

\textbf{Conclusions.} For algebras $A_{t}\left( n\right) ,$ given in
relation $\left( 2.1\right) ,$ we obtained in [Fl; 13], Theorem 2.7 and
Theorem 2.5, bounds for the level and sublevel of these algebras, namely: $%
s\left( A_{t}\left( n\right) \right) \in \lbrack n-[\frac{n}{2^{t}}],n]$ and 
$\underline{s}\left( A_{t}\left( n\right) \right) \in \lbrack n-[\frac{%
n+2^{t}-1}{2^{t}}],n],$ where $n\in \mathbb{N}-\{0\},~t\geq 2.$ From these
results, was obtained the first example of algebras obtained by the
Cayley-Dikson process of a prescribed level $n,$ in [Fl; 13], Theorem 2.9,
and a prescribed sublevel $n$ in the above Theorem 2.4. This technique
provide us an example of such algebras $A_{t}\left( n\right) $ of dimension $%
2^{t}$ and a prescribed level $n$ less than $2^{t}.$ For example, for
dimension $2^{3},$ we obtain octonion division algebras of level and
sublevel $n,$ with $n\in \{2,3,4,5,6,7,8\}.$This is the first example of
octonion algebra $\mathbb{O}$ of sublevel $6$ and sublevel$~7,$ values which
are not of the form $2^{k}$ or $2^{k}+1.$

We can't provide, using this technique, algebras obtained by the
Cayley-Dikson process of dimension $2^{t}$ \ and level greater that $2^{t}.$

This technique is as an elevator which can ascend but can not descend, since
we can't find the level and sublevel of the quaternion subalgebra $\mathbb{H}
$ of the algebra $\mathbb{O},$ in the case when $s\left( \mathbb{O}\right) ,%
\underline{s}\left( \mathbb{O}\right) \in \{6,7\}.$

In this paper, we developed another technique, which allowed us to find
better bounds for the level and sublevel of alebras obtained by the
Cayley-Dikson process of dimension $2^{t}\ $and the level $n$ greater than $%
2^{t}.$ This new technique is based on finding a new field on which the
defined algebras can have better bounds and give us help to find a positive
answer to the following question: for any positive integer $n$, how can the
existence of an algebra obtained by the Cayley-Dickson process, of dimension 
$2^{t},t\geq 4$ and level $n,$ influence the existence of a quaternion or an
octonion division algebra of level $n.$ The answer at this question can be
the key for solving the problem of the existence of quaternion and octonion
division algebras of prescribed level and sublevel.

\bigskip 
\begin{equation*}
\end{equation*}%
\begin{equation*}
\end{equation*}

\begin{equation*}
\end{equation*}

\textbf{References}

\begin{equation*}
\end{equation*}

[Br; 67] Brown, R. B., \textit{On generalized Cayley-Dickson algebras},
Pacific J. of Math.,\textbf{\ 20(3)}(1967), 415-422.\smallskip

[Da, La, Pe; 80] Dai, Z.D., Lam, T. Y., Peng, C. K., \textit{Levels in
algebra and topology,} Bull. Amer. Math. Soc., \textbf{3}(1980),845-848.%
\smallskip

[El, La; 72] Elman, R., Lam, T.Y., \textit{Pfister forms and K-theory of
fields}, Journal of Algebra \textbf{23}(1972), 181--213.

[Fl; 11] Flaut, C., \textit{Isotropy of some quadratic forms and its
applications on levels and sublevels of algebras,} J. Math. Sci. Adv. Appl., 
\textbf{12(2)}(2011), 97-117.\smallskip

[Fl; 13] Flaut, C., \textit{Levels and sublevels of algebras obtained by the
Cayley--Dickson process,} Ann. Mat. Pura Appl., \textbf{192(6})(2013),
1099-1114.

[Hoff; 95] Hoffman, D. W.,\textit{\ Isotropy of quadratic forms over the
function field of a quadric}, Math. Z, \textbf{220(3)}(1995), 461-476.

[Hoff; 08] Hoffman, D. W., \textit{Levels of quaternion algebras}, Archiv
der Mathematik, \textbf{90(5)}(2008), 401-411.\smallskip

[Ka, Me; 03] Karpenko, N.A., Merkurjev, A.S., \textit{\ Essential dimension
of quadratics}, Inventiones Mathematicae, \textbf{153}(2003),
361-372.\smallskip

[Kn; 76] Knebusch, M., \textit{Generic splitting of quadratic forms I},
Proc. London Math. Soc. \textbf{33}(1976), 65-93.\smallskip

[La, Ma; 01] Laghribi A., Mammone P., \textit{On the level of a quaternion
algebra}, Comm. Algebra, \textbf{29(4)}(2001), 1821-1828.\smallskip

[Lew; 87] Lewis, D. W., \textit{Levels and sublevels of division algebras},
Proc. Roy. Irish Acad. Sect. A, \textbf{87(1)}(1987), 103-106.\smallskip

[Lew; 89] Lewis, D. W., \textit{Levels of quaternion algebras}, Rocky
Mountain J, Math. \textbf{19}(1989), 787-792.\smallskip

[O' Sh; 06] O' Shea, J., \textit{New values for the levels and sublevels of
composition algebras}, preprint.

[O' Sh; 07] O' Shea, J., \textit{Levels and sublevels of composition algebras%
}, Indag. Mathem., \textbf{18(1)}(2007), 147-159.\smallskip

[O' Sh; 10] O' Shea, J., \textit{Bounds on the levels of composition algebras%
}, Mathematical Proceedings of the Royal Irish Academy \textbf{110A(1)}%
(2010), 21-30.\smallskip

[O' Sh; 11] O' Shea, J., \textit{Sums of squares in certain quaternion and
octonion algebras}, C.R. Acad. Sci. Paris S\'{e}r. I Math, \textbf{349}%
(2011), 239-242.\smallskip

[Pf; 65] Pfister, A., \textit{Zur Darstellung von-I als Summe von quadraten
in einem K\"{o}rper}, J. London Math. Soc. \textbf{40}(1965),
159-165.\smallskip

[Pu; 05] Pumpl\"{u}n, S., \textit{Sums of squares in octonion algebras},
Proc. Amer. Math. Soc., \textbf{133}(2005), 3143-3152.\smallskip

[Sc; 66] Schafer, R. D., \textit{An Introduction to Nonassociative Algebras,}
Academic Press, New-York, 1966.\smallskip

[Sc; 54] Schafer, R. D., \textit{On the algebras formed by the
Cayley-Dickson process,} Amer. J. Math., \textbf{76}(1954),
435-446.\smallskip

[Sch; 85] Scharlau,W., \textit{Quadratic and Hermitian Forms, }Springer
Verlag, 1985.\smallskip

[Ti, Va; 87] Tignol, J.-P., Vast, N., \textit{Representation de -1 comme
somme de carr\'{e} dans certain alg\`{e}bres de quaternions}, C.R. Acad.
Sci. Paris S\'{e}r. Math. 305, \textbf{13}(1987), 583-586.

\ 

\end{document}